\documentclass[a4paper,11pt]{article}
\usepackage[nointlimits]{amsmath}
\usepackage{amssymb,amsmath,amsthm}
\usepackage{amsfonts}
\usepackage{latexsym}
\usepackage{graphicx}
\usepackage{color}
\usepackage{layout}
\usepackage[english]{babel}
\numberwithin{equation}{section} 

\def\ds{\displaystyle}

\def\CC{\mathbb C}
\newcommand{\bin}[2]{\genfrac{(}{)}{0pt}{}{#1}{#2}}
\newcommand{\leg}[2]{\left({#1\over #2}\right)}
\newcommand{\qbin}[2]{\genfrac{[}{]}{0pt}{}{#1}{#2}_q}
\newcommand{\qqbin}[2]{\genfrac{[}{]}{0pt}{}{#1}{#2}_{1/q}}
\def\ds{\displaystyle}

\newtheorem{thm}{Theorem}[section]
\newtheorem{lem}[thm]{Lemma}
\newtheorem{cor}[thm]{Corollary}

\begin{document}

\title{\bf $q$-Analogs of some congruences\\ involving Catalan numbers}

\author{{\sc Roberto Tauraso}\\
Dipartimento di Matematica\\
Universit\`a di Roma ``Tor Vergata'', Italy\\
{\tt tauraso@mat.uniroma2.it}\\
{\tt http://www.mat.uniroma2.it/$\sim$tauraso}
}

\date{}

\maketitle

\begin{abstract}
\noindent We provide some variations on the Greene-Krammer's identity which involve $q$-Catalan numbers.
Our method reveals a curious analogy between these new identities and some congruences modulo a prime.
\end{abstract}

\makeatletter{\renewcommand*{\@makefnmark}{}
\makeatother}

\section{Introduction} 
In \cite{Gr}, John Greene proved the following conjecture made by Daan Krammer
$$
1+2\sum_{k=1}^{n-1}(-1)^{k}q^{-\bin{k}{2}}\qbin{2k-1}{k}=\left\{
\begin{array}{ll}
\ds\leg{m}{5}\sqrt{5}   & \mbox{if $5\mid n$}\\\\
\ds\leg{n}{5} & \mbox{otherwise}
\end{array}
\right.
$$
where  $q=e^{2\pi m i/n}$ with $\gcd(n,m)=1$ and $\leg{n}{p}$ is the standard Legendre symbol
(see also \cite{A2} and \cite{CS}). On the other hand if we take $q=1$, and let $n$ be a power of a prime $p$ 
then the l.h.s satisfies the following congruence which appears in \cite{ST2}
$$
1+2\sum_{k=1}^{p^a-1}(-1)^{k}\bin{2k-1}{k}=\sum_{k=0}^{p^a-1}(-1)^{k}\bin{2k}{k}\equiv\leg{p^a}{5} \pmod{p}.
$$
In this note we would like to present more examples of the same flavour
involving the $q$-Catalan numbers.

\section{Notations and properties of $q$-binomial coefficients}
\noindent The {\sl Gaussian $q$-binomial coefficient} $\qbin{n}{k}$ is defined as
$$\qbin{n}{k}=\left\{
\begin{array}{ll}
(q;q)_n(q;q)_k^{-1}(q;q)_{n-k}^{-1} &\mbox{if $0\leq k\leq n$}\\
0 &\mbox{otherwise}
\end{array}\right.$$
where $(z;q)_n=\prod_{j=0}^{m-1}(1-zq)$.
It is a polynomial in $q$ which satisfies the following relations for $0\leq k\leq n$
\begin{eqnarray}
\label{bi1}\qbin{n}{k}&=&q^{n-k}\qbin{n-1}{k-1}+\qbin{n-1}{k} \\
\label{bi2}\qbin{n}{k}&=&\qbin{n-1}{k-1}+q^{k}\qbin{n-1}{k}\\
\label{bi3}\qbin{n}{k}&=&q^{k(n-k)}\qqbin{n}{k}.
\end{eqnarray}
We define the $q$-Fibonacci polynomials (\cite{Ca}) by the recursion
$$F_n^q(t)=F_{n-1}^q(t)+q^{n-2}t F_{n-2}^q(t)$$
with initial values $F_0^q(t)=0$, $F_1^q(t)=1$.
The following identity yields an explicit formula
\begin{equation}\label{qfd}
F_n^q(t)=\sum_{k\geq 0}q^{k^2}\qbin{n-1-k}{k}t^k.
\end{equation}
There are various $q$-analogs of the Catalan numbers $C_n={1\over n+1}\bin{2n}{n}$ (see for example \cite{FH}).
We will consider the following definition
$$C_n^q={1\over [n+1]_q}\qbin{2n}{n}=\qbin{2n}{n}-q\qbin{2n}{n+1}.$$
where $[n+1]_q=(1-q)/(1-q^{n+1})$.
\noindent By \cite{CH}, $C_n^q$ is a polynomial with respect to $q$.

\section{$q$-Binomial coefficient congruences}
Let $\Phi_n(q)$ be the $n$-cyclotomic polynomial:
$$\Phi_n(q)=\prod_{\begin{array}{c}0\leq m<n\\\gcd(m,n)=1\end{array}}(q-e^{2\pi mi/n}).$$
We now deduce some properties that we will need later.

\begin{lem} For $n>1$
\begin{equation}\label{cpi} 
\Phi_n(1)=\left\{
\begin{array}{ll}
\ds p   & \mbox{if $n$ is a power of a prime $p$}\\\\
\ds 1   & \mbox{otherwise}
\end{array}
\right. .
\end{equation}
\end{lem}
\begin{proof} See for example \cite{Na} at page 160.\end{proof}
 
\begin{lem} For any positive integer $a$
\begin{equation}\label{bc1} 
\qbin{an}{k}\equiv\left\{
\begin{array}{ll}
\ds \bin{a}{k/n}   & \mbox{if $n|k$}\\\\
\ds 0   & \mbox{otherwise}
\end{array}
\right. \pmod{\Phi_n(q)}
\end{equation}
and
\begin{equation}\label{bc2}
\qbin{n+1}{k}\equiv\left\{
\begin{array}{ll}
\ds 1   & \mbox{if $k=0,1,n,n+1$}\\\\
\ds 0   & \mbox{otherwise}
\end{array}
\right. \pmod{\Phi_n(q)}.
\end{equation}
\end{lem}
\begin{proof}
By \cite{CH}, $\Phi_n(q)$ is a factor of $\qbin{m}{k}$ if and only if
$\{k/n\}>\{m/n\}$ where $\{x\}$ denote the fractional part of $x$, namely 
$\{x\}=x-\lfloor x \rfloor$. Morever, by \cite{Cl},
$$\qbin{an}{bn}\equiv\bin{a}{b}\pmod{\Phi_n(q)}.$$
\end{proof}

\begin{lem} The following congruences hold:

\noindent for $k=1,\dots,n-1$
\begin{equation}\label{bc3}
\qbin{2k-1}{k}\equiv(-1)^kq^{{3k^2-k\over 2}} \qbin{n-k}{k} \pmod{\Phi_n(q)},
\end{equation}
for $k=0,\dots,n-1$
\begin{equation}\label{bc4}
\qbin{2k}{k}\equiv(-1)^kq^{{3k^2+k\over 2}} \qbin{n-1-k}{k}\pmod{\Phi_n(q)},
\end{equation}
and
\begin{equation}\label{bc5}
\qbin{2k}{k+1}\equiv\left\{
\begin{array}{ll}
\ds (-1)^{k+1}q^{{3k^2+3k\over 2}} \qbin{n-k}{k+1}   & \mbox{if $k=0,1,\dots,n-2$}\\\\
\ds 1    & \mbox{if $k=n-1$}
\end{array}
\right. \pmod{\Phi_n(q)}.
\end{equation}
\end{lem}
\begin{proof}  Let $q=e^{2\pi m i/n}$ with $\gcd(n,m)=1$.
Since $q^k\not =1$ for $k=1,\dots,n-1$ and $q^n=1$, then 
\begin{eqnarray*}
\qbin{2k-1}{k}&=&{(1-q^{2k-1})\cdots (1-q^k)\over (1-q^{k})\cdots(1-q)}\\
&=&(-1)^kq^{{3k^2-k\over 2}}{(1-q^{n-(2k-1)})\cdots (1-q^{n-k})\over (1-q^{k})\cdots(1-q)}\\
&=&(-1)^kq^{{3k^2-k\over 2}} \qbin{n-k}{k}.
\end{eqnarray*}
Hence 
$$\qbin{2k-1}{k}-(-1)^kq^{{3k^2-k\over 2}} \qbin{n-k}{k}$$
is a the polynomial in $q$ which has at least the same roots of $\Phi_n(q)$ 
and the proof of (\ref{bc3}) is complete.
In a similar way we show the other two congruences (\ref{bc4}) and (\ref{bc5}).
\end{proof}

\section{$q$-Identities}
A fundamental result that we are going to use is the finite form of the Rogers-Ramanujan identities
(see for example \cite{A1} p. 50): for $a\in\{0,1\}$
\begin{eqnarray}\label{qid1}
F_{n+1-a}^q(q^a)=\sum_{k\geq 0}q^{k^2+ak}\qbin{n-a-k}{k}=
\sum_{j=-\infty}^{\infty}(-1)^j q^{{j(5j+1-4a)\over 2}}
\qbin{n}{\lfloor{n+2a-5j\over 2}\rfloor}
\end{eqnarray}

The next $q$-identity has been proved in \cite{EZ} with a computer proof. 
Here we show that the identity holds by using only some basic 
properties of $q$-binomial coefficients.
\begin{thm} 
\begin{equation}\label{qid2}
\sum_{k\geq 0}(-1)^k q^{\bin{k}{2}}\qbin{n-k}{k}=(-1)^n\leg{n+1}{3}q^{{1\over 3}\bin{n}{2}}
\end{equation}
\end{thm}

\begin{proof} Let the l.h.s. be $G(n)$ and let the r.h.s be $H(n)$.
It's easy to verify that $G(n)=H(n)$ for $n=0,1,2,3$. Moreover for $n\geq 1$
$$H(n+3)=-(-1)^{n+1}\leg{n+1}{3}q^{{1\over 3}\bin{n}{2}+n+1}
=-q^{n+1}H(n).$$
So it suffices to show that the same identity holds also for $G(n)$.
By (\ref{bi1}), we have that
\begin{eqnarray*}
G(n+3)&=&1-\sum_{k\geq 1}(-1)^{k-1}q^{\bin{k}{2}+n+3-2k}\qbin{n+2-k}{k-1}
+\sum_{k\geq 1}(-1)^k q^{\bin{k}{2}}\qbin{n+2-k}{k}\\
&=&G(n+2)-q^{n+1}\sum_{k\geq 1}(-1)^{k-1}q^{\bin{k-1}{2}-(k-1)}\qbin{n+1-(k-1)}{k-1}\\
&=&G(n+2)-q^{n+1}\sum_{k\geq 0}(-1)^{k}q^{\bin{k}{2}-k}\qbin{n+1-k}{k}.
\end{eqnarray*}
Moreover, by (\ref{bi2}),
\begin{eqnarray*}
\sum_{k\geq 0}(-1)^{k}q^{\bin{k}{2}-k}\qbin{n+1-k}{k}
&=& 1-\sum_{k\geq 1}(-1)^{k-1}q^{\bin{k-1}{2}-1}\qbin{n-k}{k-1}+\sum_{k\geq 1}(-1)^k q^{\bin{k}{2}}\qbin{n-k}{k}\\
&=&-q^{-1}G(n-1)+G(n).
\end{eqnarray*}
Thus, since by the induction hypothesis $G(n+2)=-q^{n}G(n-1)$ then
$$G(n+3)=G(n+2)-q^{n+1}(-q^{-1}G(n-1)+G(n))=-q^{n+1}G(n).$$
\end{proof}

The last $q$-identity seems to be new.
It is a $q$-analogue of a binomial identity
contained in \cite{ST2} and it has been conjectured for $d=0$ by
Z. W. Sun. Since the proof is rather technical we postpone it to the last section.

\begin{thm}\label{MT} For $n\geq |d|$ then
\begin{equation}\label{qid3}
\sum_{k=0}^{n-1}q^{k}\qbin{2k}{k+d}=\sum_{k=0}^{n-|d|}
q^{{1\over 3}(2(n-k)^2-(n-k)\leg{n-|d|-k}{3}-2d^2-1)}
\leg{n-d-k}{3}\qbin{2n}{k}
\end{equation}
\end{thm} 

The $p$-congruence of the next corollary has been proved in \cite{ST1} 
(see \cite{PS} for the case $a=1$).

\begin{cor}\label{ccc} Let $n\geq |d|$ then
\begin{equation}\label{qc1}
\sum_{k=0}^{n-1}q^{k}\qbin{2k}{k+d}\equiv
\leg{n-|d|}{3}q^{{3\over 2}r(r+1)+|d|(2r+1)}\pmod{\Phi_n(q)}
\end{equation}
where $r=\lfloor 2(n-|d|)/3\rfloor$.
Moreover for $a>0$ and for any prime $p$ then
$$\sum_{k=0}^{p^a-1}\bin{2k}{k+d}\equiv\leg{p^a-|d|}{3} \pmod{p}.$$
\end{cor}

\begin{proof}
For $n\geq |d|$ we use (\ref{qid3}) and, since by (\ref{bc1}) $\qbin{2n}{k}$ is $0$ modulo $\Phi_n(q)$
unless $k=0,n,2n$, then
\begin{equation*}
\sum_{k=0}^{n-1}q^{k}\qbin{2k}{k+d}\equiv
q^{{1\over 3}(2n^2-n\leg{n-|d|}{3}-2d^2-1)}\leg{n-|d|}{3}\qbin{2n}{0}\pmod{\Phi_n(q)}
\end{equation*}
and the result follows.  For the $p$-congruence, let $q=1$ and $n=p^a$ in (\ref{qc1}) and use (\ref{cpi}).
\end{proof}

\section{A dual of Greene-Krammer's identity}
By Corollary \ref{ccc}, if we take $d=0$ we have that for any prime $p$
$$1+2\sum_{k=1}^{p^a-1}\bin{2k-1}{k}=\sum_{k=0}^{p^a-1}\bin{2k}{k}\equiv\leg{p^a}{3} \pmod{p}.$$
The analogy mentioned at the beginning guided us to the following statement.

\begin{thm} Let $q=e^{2\pi m i/n}$ with $\gcd(n,m)=1$ then
$$
1+2\sum_{k=1}^{n-1}q^{k}\qbin{2k-1}{k}=\left\{
\begin{array}{ll}
\ds \leg{m}{3} i\sqrt{3}  & \mbox{if $3\mid n$}\\\\
\ds\leg{n}{3} & \mbox{otherwise}
\end{array}
\right..
$$
\end{thm}

\begin{proof} We first note that, by (\ref{bi3}),
$$q^{k}\qbin{2k-1}{k}=q^{k^2}\qqbin{2k-1}{k}
=\mbox{conj}\left(q^{-k^2}\qbin{2k-1}{k}\right)$$\\
where $\mbox{conj}(z)$ is the complex conjugate of $z\in\CC$.

\noindent Since $\Phi_n(q)=0$, by (\ref{bc3}) we have that
$$q^{-k^2}\qbin{2k-1}{k}=(-1)^k q^{\bin{k}{2}}\qbin{n-k}{k}.$$
Hence
\begin{eqnarray*}
1+2\sum_{k=1}^{n-1}q^{k}\qbin{2k-1}{k}
&=&\mbox{conj}\left(1+2\sum_{k=1}^{n-1}q^{-k^2}\qbin{2k-1}{k}\right)\\
&=&\mbox{conj}\left(-1+2\sum_{k\geq 0} (-1)^k q^{\bin{k}{2}}\qbin{n-k}{k}\right)\\
&=&-1+2(-1)^n\leg{n+1}{3}q^{-{1\over 3}\bin{n}{2}}
\end{eqnarray*}
and the result is easily deduced.
\end{proof}

\section{$q$-Catalan congruences}

\begin{thm} For $n>0$
\begin{equation}\label{C3}
\sum_{k=0}^{n-1}q^{k}C_k^q\equiv\left\{
\begin{array}{lll}
\ds q^{\lfloor n/3\rfloor}  & \mbox{if $n=0,1$}&\pmod{3}\\\\
\ds-1-q^{(2n-1)/3} & \mbox{if $n=2$}&\pmod{3}
\end{array}
\right.\pmod{\Phi_n(q)}.
\end{equation} 
\end{thm}
\begin{proof}
Since
$$C_n^q=\qbin{2n}{n}-q\qbin{2n}{n+1}$$
by using (\ref{qc1}) for $d=0$ and for $d=1$ we obtain
\begin{equation*}
\sum_{k=0}^{n-1}q^{k}C_k^q\equiv
q^{{1\over 3}(2n^2-n\leg{n}{3}-1)}\leg{n}{3}
-q^{{1\over 3}(2n^2-n\leg{n-1}{3})}\leg{n-1}{3}\pmod{\Phi_n(q)}.
\end{equation*}
Finally we proceed by cases on $n$ modulo $3$ and the proof is complete.
\end{proof}

The $p$-congruence of the next corollary has been proved in \cite{ST1} 
(see \cite{PS} for the case $a=1$).

\begin{cor} Let $q=e^{2\pi m i/n}$ with $\gcd(n,m)=1$.  
If $3\mid n$ then
\begin{equation*}
\sum_{k=0}^{n-1}q^{k}C_k^q={1\over 2}\left(i\sqrt{3}\leg{m}{3}-1\right) 
\end{equation*}
Moreover, for any prime $p$ and for $a>0$
$$\sum_{k=0}^{p^a-1}C_k\equiv{1\over 2}\left(3\leg{p^a}{3}-1\right) \pmod{p}.$$
\end{cor}
\begin{proof} If $3\mid n$ then
$$\sum_{k=0}^{n-1}q^kC_k^q=q^{n/3}=e^{2\pi m i/3}={1\over 2}\left(i\sqrt{3}\leg{m}{3}-1\right).$$
As regards the $p$-congruence, let $q=1$ and $n=p^a$ in (\ref{C3}) and use (\ref{cpi}).
\end{proof}

\begin{thm} For $n>0$
\begin{equation}\label{C5}
\sum_{k=0}^{n-1}(-1)^k q^{-\bin{k}{2}}C_k^q\equiv
 F_n^q(q)+F_{n+2}^q(1)-2\pmod{\Phi_n(q)}.
\end{equation} 
and 
$$F_n^q(q)+F_{n+2}^q(1)\equiv \left\{
\begin{array}{lll}
\ds (-1)^{r(n)}\left(q^{r(n)(n-1)\over 2}+q^{r(n)(n+1)\over 2}\right) & \mbox{if $n\equiv 0,2,3$} &\pmod{5}\\\\
\ds (-1)^{r(n)}\left(q^{r(n)(n-2)\over 2}+q^{r(n)n\over 2}+q^{r(n)(n+2)\over 2}\right)   & \mbox{if $n\equiv 1,4$}&\pmod{5}
\end{array}\right.$$
where $r(n)=\mbox{round}(n/5)=\lfloor n/5+{1\over 2}\rfloor$.
\end{thm}
\begin{proof} By  (\ref{bc4}) and (\ref{bc5})
\begin{eqnarray*}
(-1)^k q^{-\bin{k}{2}}C_k^q
&\equiv&(-1)^k q^{{-k^2+k\over 2}}\qbin{2k}{k}
-(-1)^k q^{{-k^2+k+2\over 2}}\qbin{2k}{k+1}\\
&\equiv&q^{k^2+k}\qbin{n-1-k}{k}
+q^{(k+1)^2}\qbin{n-1-(k+1)}{k+1}\\
&&-[k=n-1] \pmod{\Phi_n(q)}.
\end{eqnarray*}
Hence by applying (\ref{qid1}) both for  $a=0$ and $a=1$ we get
\begin{eqnarray*}
\sum_{k=0}^{n-1}(-1)^k q^{-\bin{k}{2}}C_k^q
&\equiv&\sum_{k\geq 0}^{n-1}q^{k^2+k}\qbin{n-1-k}{k}+
\sum_{k\geq 0}^{n-1}q^{(k+1)^2}\qbin{n-k}{k+1}-1\\
&\equiv&
\sum_{k\geq 0}q^{k^2+k}\qbin{n-1-k}{k}+
\sum_{k\geq 0}^{n-1}q^{k^2}\qbin{n+1-k}{k}-2\\
&\equiv&
F_n^q(q)+F_{n+2}^q(1)-2 \pmod{\Phi_n(q)}.
\end{eqnarray*}

By (\ref{qid1}) for $a=1$ we find
$$F_n^q(q)=\sum_{j=-\infty}^{\infty}(-1)^j
q^{{j(5j-3)\over 2}}\qbin{n}{\lfloor{n+2-5j\over 2}\rfloor}.$$
Thus, by proceeding by cases on $n$ modulo $5$ we obtain
\begin{equation*}\label{qfib1}
F_n^q(q)\equiv\left\{
\begin{array}{lll}
\ds (-1)^{r(n)}q^{(n+2)(n-1)\over 10} &\mbox{if $n=1,3$} &\pmod{5}\\\\
\ds (-1)^{r(n)}q^{(n+1)(n-2)\over 10} &\mbox{if $n=2,4$} &\pmod{5}\\\\
\ds 0 & \mbox{if $n=0$} &\pmod{5}
\end{array}
\right.\pmod{\Phi_n(q)}
\end{equation*}

Similarly by (\ref{qid1}) for $a=0$
$$
F_{n+2}^q(1)=\sum_{j=-\infty}^{\infty}(-1)^j q^{{j(5j+1)\over 2}}\qbin{n+1}{\lfloor{n+1-5j\over 2}\rfloor}
$$
and 
\begin{equation*}
F_{n+2}^q(1)\equiv\left\{
\begin{array}{lll}
\ds (-1)^{r(n)}\left(q^{r(n)(n+1)\over 2}+q^{r(n)(n-1)\over 2}\right) &\mbox{if $n=0$} &\pmod{5}\\\\
\ds (-1)^{r(n)}\left(q^{r(n)n\over 2}+q^{r(n)(n-2)\over 2}\right) &\mbox{if $n=1$} &\pmod{5}\\\\
\ds (-1)^{r(n)}q^{r(n)(n-1)\over 2}                               &\mbox{if $n=2$} &\pmod{5}\\\\
\ds (-1)^{r(n)}q^{r(n)(n+1)\over 2}                               &\mbox{if $n=3$} &\pmod{5}\\\\
\ds (-1)^{r(n)}\left(q^{r(n)n\over 2}+q^{r(n)(n+2)\over 2}\right) &\mbox{if $n=4$} &\pmod{5}
\end{array}
\right.\pmod{\Phi_n(q)}
\end{equation*}
\end{proof}
The $p$-congruence of the next corollary has been proved in \cite{ST2}. 

\begin{cor} Let $q=e^{2\pi m i/n}$ with $\gcd(n,m)=1$.
If $5\mid n$ then
\begin{equation*}
\sum_{k=0}^{n-1}(-1)^k q^{-\bin{k}{2}}C_k^q={1\over 2}\left(\sqrt{5}\leg{m}{5}-3\right).
\end{equation*}
Moreover, for $a>0$ and for any prime $p$
$$\sum_{k=0}^{p^a-1}(-1)^kC_k\equiv{1\over 2}\left(5\leg{p^a}{5}-3\right) \pmod{p}.$$
\end{cor}
\begin{proof}
If $5\mid n$ then by (\ref{C5}):
\begin{eqnarray*}
\sum_{k=0}^{n-1}(-1)^k q^{-\bin{k}{2}}C_k^q&=&
-2+(-1)^{n/5}\left(q^{{n(n-1)\over 10}}+q^{{n(n+1)\over 10}}\right)\\
&=&-2+(-1)^{n(m+1)/5}2\,\mbox{Re}(e^{\pi i m/5})\\
&=&{1\over 2}\left(\sqrt{5}\leg{m}{5}-3\right).
\end{eqnarray*}
The $p$-congruence follows by letting $q=1$, $n=p^a$ in (\ref{C5})
and by noting that
$$(-1)^{r(p^a)}=\leg{p^a}{5}\qquad\mbox{for $p\not=5$}$$
then use (\ref{cpi}).
\end{proof}

\section{Proof of Theorem \ref{MT}}
Let
$$S(n,d)=\sum_{k=0}^{n-1}q^k\qbin{2k}{k+d}.$$
The finite sum $S(n,d)$ has some interesting properties. 

\noindent The first one is that 
$$S(n,d)=S(n,-d)$$ 
which follows immediately
from $\qbin{2k}{k+d}=\qbin{2k}{k-d}$. 
The second one is less trivial.

\begin{lem} Let $n\geq |d|$ then
\begin{equation}\label{3shift}
S(n,d)-q^{4d+6} S(n,d+3)=q^d {[2d+3]_q \over [2n+1]_q}\qbin{2n+1}{n+d+2}-[d=-1]q^{-1}+[d=-2]q^{-3}.
\end{equation}
\end{lem}

\begin{proof}
We first consider the case when $d\geq 0$ and we prove (\ref{3shift}) 
by induction on $n$. 

\noindent For $n=1$ it holds. Now we are going to prove that for $n\geq 1$
$$S(n+1,d)-q^{4d+6}S(n+1,d+3)=q^d {[2d+3]_q \over [2n+3]_q}\qbin{2n+1}{n+d+3}.$$
Since the l.h.s. is equal to
$$S(n,d)+q^n\qbin{2n}{n+d}-q^{4d+6}\left(S(n,d+3)+q^n\qbin{2n}{n+d+3}\right),$$
by the induction hypothesis, it suffices to show that
$$q^{n-d}\qbin{2n}{n+d}-q^{n+3d+6}\qbin{2n}{n+d+3}=
{[2d+3]_q\over [2n+3]_q}\qbin{2n+1}{n+d+3}-{[2d+3]_q\over [2n+1]_q}\qbin{2n+1}{n+d+2}$$
which holds.
Since 
$$q^{k+1}\qbin{2k}{k+1}-q^{k+3}\qbin{2k}{k+2}=C^q_{k+1}-C^q_{k}$$
then 
$$qS(n,1)-q^3S(n,2))=C^q_n-1$$ 
and (\ref{3shift}) holds for $d=-1$
$$S(n,-1)-q^{2} S(n,2)=S(n,1)-q^{2} S(n,2)=q^{-1}C_n^q-q^{-1},$$
and for $d=-2$
$$S(n,-2)-q^{-2} S(n,1)=S(n,2)-q^{-2} S(n,1)=-q^{-3}C_n^q+q^{-3}.$$
If $d\leq -3$, by letting $d'=-d-3\geq 0$ we get
\begin{eqnarray*}
S(n,d)-q^{4d+6}S(n,d+3)&=&S(n,-d)-q^{4d+6}S(n,-d-3)\\
&=&-q^{-4d'-6}(S(n,d') -q^{4d'+6}S(n,d'+3))\\
&=&-q^{-3d'-6}{[2d'+3]_q\over [2n+1]_q}\qbin{2n+1}{n+d'+2}\\
&=&q^{d}{[2d+3]_q\over [2n+1]_q}\qbin{2n+1}{n+d+2}.
\end{eqnarray*}
\end{proof}

Finally we are ready to prove Theorem \ref{MT}.
\begin{proof}
The $q$-identity $(\ref{qid3})$ is equivalent to $S(n,d)=T(n,d)$ where
$$T(n,d)=
\sum_{k\geq 0} q^{6k^2+(3+|d|)k+|d|}\qbin{2n}{n+3k+|d|+1}
-\sum_{k\geq 1} q^{6k^2+(3+|d|)k+|d|}\qbin{2n}{n+3k+|d|+1}.$$
By the previous lemma it suffices to verify it for $d=0,1$
(remember that $S(n,1)=S(n,-1)$). 
Since the proof for $d=1$ is quite similar, we will consider only the case for $d=0$.
So
$$T(n,0)=\sum_{k=0}^{\infty}s(n,k,3,1)-\sum_{k=1}^{\infty}s(n,k,-3,-1)$$
where 
$$s(n,k,a,b)=q^{6k^2+ak}\qbin{2n}{n+3k+b}.$$

\noindent By using the Maple package $q$-Zeilberger, we verified that 
$s(n,k,a,b)$ solves the following recurrence
$$\sum_{j=0}^4 a_j(n,a,b) s(n+j,k,a,b)=g(n,k+1,a,b)-g(n,k,a,b)$$
where $g(n,k,a,b)=r(n,k,a,b)s(n,k,a,b)$,
\begin{eqnarray*}
a_0(n,a,b)&=&(1-q^{2n+1})(1-q^{2n+2})q^6\\
a_1(n,a,b)&=&-(q^{4n+7+a-4b}+q^{4n+7-a+4b}-q^{2n+4}-q^{2n+3}+q^3+q^2+q+1)q^3\\
a_2(n,a,b)&=&(q^{4n+10}+q^{2n+7}+q^{2n+6}+q^{4}+q^3+2q^2+q+1)q\\
a_3(n,a,b)&=&-(q^{2n+6}+q^{2}+1)(q+1)\\
a_4(n,a,b)&=&1
\end{eqnarray*}
and $r(n,k,a,b)$ is the rational function certificate
$${h(n,k,a,b)(1-q^{2n+1})(1-q^{2n+2})q^{4n-12k+10-a-4b} \over
(1-q^{n-3k-b+3})(1-q^{n-3k-b+4})(1-q^{n-3k-b+1})(1-q^{n+3k+b+1})}$$
with
\begin{eqnarray*} 
h(n,k,a,b)&=&
+q^{3n+9k+3b+a+6}-q^{3n+3k+5b+9}-q^{2n+6k+2b+a+7}+q^{2n+6k+6b+7}\\
&&-q^{2n+6k+2b+a+6}+q^{2n+6k+6b+6}+q^{n+3k+b+a+7}+q^{n+3k+b+a+6}\\
&&-q^{2n+6k+2b+a+5}+q^{2n+6k+5+6b}-q^{n+9k+7b+4}+q^{n+3k+b+a+5}\\
&&-q^{n+9k+7b+3}-q^{n+9k+7b+2}+q^{12k+8b}-q^{a+6}.
\end{eqnarray*}
\noindent Hence, since $g(n,k,a,b)=r(n,k,a,b)s(n,k,a,b)=0$ when $|3k+b|>n$,
$$\sum_{j=0}^4 a_j(n,a,b)\sum_{k=k_0}^{\infty}s(n+j,k,a,b)=
\sum_{k=k_0}^{\infty}\left( g(n,k+1,a,b)-g(n,k,a,b)\right)=-g(n,k_0,a,b).$$
Since $a_j(n,3,1)=a_j(n,-3,-1)$ then
\begin{eqnarray*}
\sum_{j=0}^4 a_j(n,3,1) T(n+j,0)&=&-g(n,0,3,1)+g(n,1,-3,-1)\\
&=&-r(n,0,3,1)\qbin{2n}{n+1}+r(n,1,-3,-1)q^{3}\qbin{2n}{n+2}.
\end{eqnarray*}
The identity $S(n,0)=T(n,0)$ holds for $n=1,2,3,4$ by direct verification. 

\noindent By induction, it holds for $n>4$ as soon as for $n>1$
$$\sum_{j=0}^4 a_j(n,3,1) S(n+j,0)
=-r(n,0,3,1)\qbin{2n}{n+1}+r(n,1,-3,-1)q^{3}\qbin{2n}{n+2}.$$
Let $c_i(n,3,1)=\sum_{j=i}^4 a_j(n,3,1)$ for $i=0,1,2,3,4$. 
Since $c_0(n,3,1)=0$, the r.h.s. can be simplified, and it suffices to check that
$$\sum_{i=0}^3 c_{i+1}(n,3,1)q^{n+i}\qbin{2(n+i)}{n+i}=-r(n,0,3,1)\qbin{2n}{n+1}+r(n,1,-3,-1)q^{3}\qbin{2n}{n+2}.$$
which holds.
\end{proof}


\end{document}